\documentclass[]{interact}

\usepackage[authoryear]{natbib}
\bibpunct{(}{)}{,}{a}{}{,}

\usepackage[english]{babel}
\usepackage[utf8]{inputenc}
\usepackage{graphicx}
\usepackage{enumitem}
\usepackage{hyperref}
\usepackage{graphicx}
\usepackage{float}
\usepackage{subcaption}
\usepackage{amsmath}
\usepackage{amsthm}
\usepackage{amssymb}
\usepackage{mathtools}
\usepackage{dsfont}

\theoremstyle{plain}
\newtheorem{theorem}{Theorem}[section]
\newtheorem{lemma}[theorem]{Lemma}
\newtheorem{corollary}[theorem]{Corollary}
\newtheorem{assumption}[theorem]{Assumption}

\theoremstyle{definition}
\newtheorem{definition}[theorem]{Definition}

\theoremstyle{remark}
\newtheorem{remark}{Remark}

\DeclareMathOperator*{\argmin}{arg\,min}
\DeclareMathOperator*{\tr}{\text{tr}}

\begin{document}

\title{Optimal Guaranteed Cost Control of Discrete-Time Linear Systems subject to Structured Uncertainties}
\author{
\name{Carlos M. Massera\textsuperscript{a}, Marco H. Terra\textsuperscript{b} and Denis F. Wolf\textsuperscript{a}}
\affil{\textsuperscript{a}Institute of Mathematics and Computer Science, University of São Paulo, São Carlos, Brazil; \textsuperscript{b}São Carlos School of Engineering, University of São Paulo, São Carlos, Brazil}
}
\thanks{Corresponding author: Carlos M. Masssera. Email: carlos.magalhaes.filho@usp.br}
\maketitle

\begin{abstract}
In this paper, we propose state- and static output-feedback generalized guaranteed cost control (GCC) approaches for discrete-time linear systems subject to norm-bounded structured parametric uncertainties. 
This method enables the convex synthesis for a more general class of systems, where uncertainties are uncorrelated block diagonal, and no feed-through uncertainty is multiplicative with control input ones.
It also provides necessary and sufficient conditions for state-feedback and sufficient conditions for static output-feedback. 
We also present a comparative study among the proposed controllers, standard Linear Quadratic Regulator, and Guaranteed Cost Controller found in the literature.
\end{abstract}

\section{Introduction}

The Linear Quadratic Regulator (LQR) \citep{kalman1960contributions} is an optimal control approach which aims to drive states of a linear system to its origin through the minimization of a quadratic cost functional. Safonov and Athans \citep{safonov1976gain} have shown that such controller has $6 dB$ gain margin and $60^\circ$ phase margin if the system dynamics is assumed to be known. However, when it is subject to uncertainties, optimal closed-loop performance cannot be obtained since robustness properties are no longer guaranteed \citep{douglas1991linear}.

To address this problem, several authors in the last decades have proposed different approaches to solve it. Petersen \citep{petersen1987stabilization}, for instance, proposed a Guaranteed Cost Controller (GCC) for a continuous-time linear system subject to $l2$-norm-bounded uncertainties based on the solution of Algebraic Riccati Equations (ARE). Such method synthesizes a controller which presents robust asymptotic stability and provides a guaranteed upper bound to a quadratic cost function used as performance measurement. Xie and Soh \citep{xie1993control} applied such approach to discrete-time systems, which Petersen \textit{et al.} \citep{petersen1998optimal} extended to a more general class of systems through the use of Linear Matrix Inequality (LMI). Previous studies applied this class of synthesis to several types of systems and controllers, including static and dynamic output feedback controllers. Some examples can be seen in: \citep{garcia2003robust, moheimani1996optimal, yu1999lmi, chen2003delay,boukas1998stochastic, yang2016guaranteed, qiu2016output, wang2016robust}.

A challenging aspect of most robust synthesis problems is its inherent non-convexity, even when no uncertainties are considered. This topic has attracted significant research interest in the past two decades. 
Most of the approaches proposed for such class of problems are based on dilated variables, first proposed by Oliveira et al. \citep{oliveira1999new} for discrete time systems and later extended by Apkarian et al. \citep{apkarian2001continuous} for continuous time systems. 
Some examples of applications of such method are: Output feedback synthesis (i.e., \citep{lee2006sufficient}, \citep{xie1996output}, \citep{dong2007static}, \citep{peaucelle2014lmi}, \citep{arzelier2003robust}); systems subject to delays (i.e., \citep{he2004parameter}, \citep{li2016delay}); systems subject to Markovian jumps (i.e. \citep{de2003mode}, \citep{morais2016reduced}, \citep{shen2015new}); constrained optimal controllers (i.e. \citep{cuzzola2002improved}, \citep{xia2008robust}, \citep{lu2013constrained}); and many others. 
The use of such dilation variables is conservative, thus it only provides sufficient conditions for the controller existence. To partially overcome this conservativeness, several authors have proposed the use of iterative methods to locally solve the resulting non-convex Bilinear Matrix Inequality (BMI) problem (i.e. \citep{sadeghzadeh2016robust}, \citep{moreira2011robust}, \citep{sadeghzadeh2015fixed}, \citep{peaucelle2001efficient}, \citep{donkers2017decentralised}, \citep{sadeghzadeh2014fixed}). However, iterative methods also only provide sufficient conditions, due to its local optimality guarantees.

Another conservative aspect of robust synthesis and analysis methods is the assumption of an arbitrary (unstructured) norm bounded set since uncertainties are often known to be structured. Several studies investigated analysis and synthesis methods for linear systems under structured uncertainties. Shamma \citep{shamma1994robust}, Packard and Doyle \citep{packard1993complex}, and Graham \textit{et al.} \citep{graham2006linear} addressed the problem of robust stability analysis. Sznaier \textit{et al.} \citep{sznaier2002convex}, and Feron \citep{feron1997analysis} investigated robust $H_2$ analysis. While Wu and Lu \citep{wu2004convexified} addressed $ H_\infty $ synthesis. However, according to Paganini \citep{Paganini2013}, the synthesis of a controller for systems subject to structured uncertainties results in non-convex optimizations.

In this paper, we propose a generalized Guaranteed Cost Control synthesis problem for discrete-time linear systems subject to norm-bounded structured parametric uncertainties for state-feedback and static output-feedback controllers. We consider a correlated block diagonal uncertainty structure, which is more general than previous studies (see for instance Szanaier \textit{et al.} \citep{sznaier2002convex} and Graham \textit{et al.} \citep{graham2006linear}). We also address the non-convexity issue discussed by Paganini through the use of the generalized S-Procedure, as proposed by Iwasaki \textit{et al.} \citep{iwasaki2000generalized}, which is proved to be both necessary and sufficient for the considered structure. Meanwhile, we address the non-convexity of the output feedback synthesis problem conservatively through the dilated variable approach, based on the Reciprocal Projection Lemma \citep{apkarian2001continuous}.

This paper has the following organization: Section \ref{sec_problem_statement} presents the problem statement; Section \ref{sec_preliminary_results} discusses preliminary results; Section \ref{sec_gcc} derives the proposed Guaranteed Cost Controller; Section \ref{sec_example} provides a numerical example; finally, Section \ref{sec_conclusion} addresses the final remarks.

The notation used throughout this paper is standard. $ I_p $ is a identity matrix of size $ p \times p $, $ x^T M x = x^T M (\bullet) $ for brevity, $ A \otimes B $ is the Kronecker product of matrices $A$ and $B$, and $ \tr(A) $ is the trace of matrix $ A $.

\section{Problem statement}
\label{sec_problem_statement}

In this section, we present a Guaranteed Cost Controller definition with its required assumptions. We consider discrete-time linear system subject to parametric uncertainties of the form
\begin{equation}
\begin{aligned}
x_{k+1} &= (A + \delta A(\Delta_k)) x_k + (B^u + \delta B^u(\Delta_k)) u_k\\
y_{k} &= (C_y + \delta C_y(\Delta_k)) x_k,
\end{aligned}
\label{eq_sys_model}
\end{equation}%
where $ x_k \in \Re^{n_x} $ is the system state, $ u_k \in \Re^{n_u} $ is the control input, $ A \in \Re^{n_x \times n_x} $ is the state matrix, $ B^u \in \Re^{n_x \times n_u} $ is the input matrix, $ C_y \in \Re^{n_y \times n_x} $ is the output state matrix, and $ \delta A(\Delta_k) $, $ \delta B^u(\Delta_k) $ and $ \delta C_y(\Delta_k) $ are, respectively, the state, input and output multiplicative uncertainty matrices, such that%
\begin{equation}
\begin{bmatrix}
\delta A(\Delta_k) & \delta B^u(\Delta_k) \\
\delta C_y(\Delta_k) & 0 \\
\end{bmatrix}
= 
\begin{bmatrix}
B^w \\ D_y^w
\end{bmatrix}
\bar{\Delta}_k
\begin{bmatrix}
C_z & D_z^u
\end{bmatrix},
\label{eq_sys_uncertainty_matrix}
\end{equation}%
where $ \bar{\Delta}_k = \Delta_k (I_{n_q} - D_z^w \Delta_k)^{-1} $, $ \Delta_k \in \mathbb{D}^u = \{ \Delta \mid \Delta \in \Re^{n_p \times n_q}, || \Delta ||_2 \le 1 \} $, $ B^w \in \Re^{n_x \times n_p} $, $ D_y^w \in \Re^{n_y \times n_p} $, $ D_z^w \in \Re^{n_q \times n_p} $, $ C_z \in \Re^{n_q \times n_x} $ and $ D_z^u \in \Re^{n_q \times n_u} $.

In the scope of this paper, we are interested in investigating the controller synthesis when the structure of $ \Delta_k $ is known. Therefore, we define $ \mathbb{D} \subseteq \mathbb{D}^u $ to be the set of uncertainty matrices such that%
\begin{multline}
\mathbb{D} = \{ \Delta \mid \forall i \in [1, s]: \Delta_i \in \Re^{n_{pi} \times n_{qi}}, || \Delta_i ||_2 \le 1, \\ \Delta = \text{diag}(I_{n_{r1}} \otimes \Delta_1, I_{n_{r2}} \otimes \Delta_2, \ldots, I_{n_{rs}} \otimes \Delta_s) \},
\label{eq_uncertainty_structure}
\end{multline}%
where $ n_p = \underset{i = 1}{\overset{s}{\sum}} n_{ri} n_{pi} $ and $ n_q = \underset{i = 1}{\overset{s}{\sum}} n_{ri} n_{qi} $.

The uncertainty structure defined in \eqref{eq_uncertainty_structure} is more general than ones considered in previous studies. For example, the structure investigated by Sznaier \textit{et al.} \cite{sznaier2002convex} and Graham \textit{et al.} \cite{graham2006linear} is square block diagonal ($ n_{ri} = 1 $ and $ n_{pi} = n_{qi} $).

Associated with the uncertain linear system \eqref{eq_sys_model} is the performance measurement cost function%
\begin{equation}
J(x_0) = \underset{n \rightarrow \infty}{\lim} \underset{k = 0}{\overset{n}{\sum}} x_k^T Q x_k + 2 x_k^T N u_k + u_k^T R u_k,
\label{eq_cost_function}
\end{equation}%
where $ Q \in \Re^{n_x \times n_x} $, $ N \in \Re^{n_x \times n_u} $, $ R \in \Re^{n_u \times n_u} $, $Q \succeq 0 $, $ R \succ 0 $, and%
\begin{equation}
\begin{bmatrix}
Q & N\\
N^T & R
\end{bmatrix} = \begin{bmatrix}
C_c^T \\ D_c^{uT}
\end{bmatrix}
\begin{bmatrix}
C_c & D_c^u
\end{bmatrix}
\succeq 0,
\end{equation}%
and $ C_c \in \Re^{n_c \times n_x} $ and $ D_c^u \in \Re^{n_c \times n_u} $ denote the factorization of the cost function.

We are now ready to define the proposed Guaranteed Cost Controller.

\begin{definition}
A static output feedback controller $u_k = -Ky_k$ is said to be a stabilizing Guaranteed Cost Controller for the uncertain system (\ref{eq_sys_model}) if there exists a symmetric matrix $P \succ 0 $ that upper-bounds the cost functional \eqref{eq_cost_function} for

\begin{multline}
[A + \delta A + (B^u + \delta B^u(\Delta_k) K (C_y + \delta C_y(\Delta_k)]^T P[\bullet] - P + \\ Q + N K (C_y + \delta C_y(\Delta_k)) + (C_y + \delta C_y(\Delta_k))^T K^T N^T + \\ + [ K (C_y + \delta C_y(\Delta_k))]^T R[\bullet] \le 0,
\end{multline}
such that $ J(x_0) \le J^*(x_0) = x_0^T P x_0 $ for the closed loop system%
\begin{equation}
x_{k+1} = \left[ A + \delta A(\Delta_k) - \left( B^u + \delta B^u(\Delta_k) \right) K \left( C_y + \delta C_y(\Delta_k) \right) \right] x_k,
\label{eq_sys_model_cl_delta}
\end{equation}%
and all admissible uncertainties $ \Delta_k \in \mathbb{D} $.
\label{def_gcc}
\end{definition}

\begin{definition} [Optimal Guaranteed Cost Control]
A static output-feedback controller $ u_k = -K y_k $ is said to be optimal if it satisfies Definition \ref{def_gcc} and $ J^*(x_0) $ is minimal.
\label{def_gcc_optimal}
\end{definition}

\begin{remark}
The static output-feedback GCC from Definition \ref{def_gcc} reduces to the state-feedback case if $ C_y = I_{n_x} $, $ D_y^w = 0 $, and $ \delta C_y(\Delta_k) = 0 $ for all $ \Delta_k \in \mathbb{D} $.
\end{remark}

The GCC problem from Definitions \ref{def_gcc} and \ref{def_gcc_optimal} is well posed under the following assumptions:

\begin{assumption}
For all $ \Delta_k \in \mathbb{D} $, the pair ($A + \delta A(\Delta_k)$, $B + \delta B(\Delta_k)$) is stabilizable.
\label{ass_stabilizable}
\end{assumption}

\begin{assumption}
The uncertainty matrices \eqref{eq_sys_uncertainty_matrix} are well-posed. Therefore, for all $ \Delta_k \in \mathbb{D} $, $ D_z^w \Delta_k \prec \mathbf{I_{n_q}}$, or equivalently, $ ||D_z^w||_2 < 1 $.
\label{ass_wellposed}
\end{assumption}

\begin{assumption}
For all $ \Delta_k \in \mathbb{D} $, the pair ($A + \delta A(\Delta_k)$, $Q^{\frac{1}{2}}$) is observable.
\label{ass_observable}
\end{assumption}

\begin{assumption}
For all $ \Delta_k \in \mathbb{D} $, $ D_y^w \bar{\Delta}_k D_z^u = 0 $. Therefore, the feed-through uncertainty is zero.
\label{ass_no_feedthrough}
\end{assumption}

\section{Preliminary results}
\label{sec_preliminary_results}

This section presents previous results and preliminary work to formulate a convex condition for the Guaranteed Cost Controller. Such results include the generalized S-Procedure, the Reciprocal Projection Lemma, a property of the Kronecker Product, and an alternative representation of the system presented in Section \ref{sec_problem_statement}.

The uncertain linear system \eqref{eq_sys_model} is equivalent to%
\begin{equation}
\begin{aligned}
x_{k+1} &= A x_k + B^w w_k + B^u u_k\\
y_k &= C_y x_k + D_y^w w_k\\
z_k &= C_z x_k + D_z^w w_k + D_z^u u_k,
\end{aligned}
\label{eq_sys_model_2}
\end{equation}%
subject to the feedback disturbance $ w_k = \Delta_k z_k $, where $ w_k \in \Re^{n_p} $ and $ z_k \in \Re^{n_q} $. Similarly, the closed loop system with a static output-feedback controller $ u_k = - K y_k $ is equivalent to%
\begin{equation}
\begin{aligned}
x_{k+1} &= (A - B^u K C_y) x_k + (B^w - B^u K D_y^w) w_k\\
z_k &= (C_z - D_z^u K C_y) x_k + (D_z^w - D_z^u K D_y^w) w_k,
\end{aligned}
\label{eq_sys_model_cl}
\end{equation}%
whose matrices will be denoted as $ \bar{A} = A - B^u K C_y $, $ \bar{B}^w = B^w - B^u K D_y^w$, $\bar{C}_z = C_z - D_z^u K C_y$, and $ \bar{D}_z^w = D_z^w - D_z^u K D_y^w $ for brevity. From \eqref{eq_sys_model_cl}, the closed loop uncertain system \eqref{eq_sys_model_cl_delta} is equivalent to%
\begin{equation}
x_{k+1} = \left( \bar{A} - \bar{B}^w \bar{\Delta}_k \bar{C}_z \right) x_k,
\label{eq_sys_model_cl_delta_2}
\end{equation}%
and the closed-loop cost function \eqref{eq_cost_function} can be rewritten as %
\begin{equation}
J(x_0) = \underset{k = 0}{\overset{\infty}{\sum}} x_k^T \left( \bar{C}_c + \bar{D}_c^w \bar{\Delta}_k \bar{C}_z \right)^T (\bullet) x_k,
\label{eq_cost_function_cl}
\end{equation}%
where $ \bar{C}_c = C_c - D_c K C_y $, and $ \bar{D}_c^w = - D_c K D_y^w $. Therefore, based on \eqref{eq_sys_model_cl_delta_2} and \eqref{eq_cost_function_cl}, the optimal GCC synthesis problem from Definition \ref{def_gcc_optimal} is equivalent to%
\begin{equation}
\begin{aligned}
& J^*(x_0) \ge \underset{\Delta_0, \ldots, \Delta_\infty}{\max} \underset{K}{\min} \;\; J(x_0)\\
& s.t. \;\; x_{k+1} = \left( \bar{A} - \bar{B}^w \bar{\Delta}_k \bar{C}_z \right) x_k.
\end{aligned}
\label{eq_gcc_opt_problem}
\end{equation}

\subsection{Generalized S-Procedure}

Iwasaki \textit{et al.} \citep{iwasaki2000generalized} described a generalized S-Procedure which provides both necessary and sufficient conditions for particular sets of matrices, entitled lossless sets. We now present the lossless set definition and the generalized S-Procedure theorem for completeness.

\begin{definition} [Lossless Set]
A subset $\mathbb{S}$ of $n \times n$ symmetric real matrices is said to be lossless if it has the following properties:%
\begin{enumerate}[label=(\roman*)]
\item $\mathbb{S}$ is convex;
\item $S \in \mathbb{S} \Rightarrow \forall \tau > 0: \tau S \in \mathbb{S}$;
\item For each nonzero matrix $H \in \Re^{n \times n}, H \succeq 0 $ such that $ \forall S \in \mathbb{S} : \tr(SH) \le 0 $, there exist vectors $ \lambda_i \in \Re^{n} $ such that%
\begin{equation}
H = \underset{i = 0}{\overset{\text{rank}(H)}{\sum}} \lambda_i \lambda_i^*,\; \forall S \in \mathbb{S} : \lambda_i^* S \lambda_i \le 0.
\end{equation}
\end{enumerate}
\label{def_lossless_set}
\end{definition}

\begin{theorem} [Generalized S-Procedure]
Let $ \Theta \in \Re^{n \times n} $ be a symmetric matrix and $\mathbb{S}$ be a lossless set. Then, the following statements are equivalent:%
\begin{enumerate}[label=(\roman*)]
\item $\forall x \in \mathbb{X} = \{x \mid x \in \Re^n, x \not = 0,x^T S x \le 0, \forall S \in \mathbb{S} \}: x^T \Theta x \le 0$;
\item $ \exists S \in \mathbb{S} : \Theta - S \preceq 0 $.
\end{enumerate}

\begin{proof}
See Theorem 1 of \citep{iwasaki2000generalized}.
\end{proof}
\label{the_generalized_s_procedure}
\end{theorem}

\subsection{Reciprocal projection lemma}

Apkarian \textit{et al.} \citep{apkarian2001continuous}, based on Oliveira et. al. \citep{de1999new}, proposed the Reciprocal Projection Lemma. Such Lemma enables several previously non-convex conditions present in synthesis problems to be posed as convex conditions through the use of additional slack variables. Based on this result less conservative methods have been developed, particularly in the field of multi-objective synthesis and analysis (i.e., \citep{korouglu2014new}, \citep{adegas2013new}). We now present the Projection and the Reciprocal Projection Lemmas.

\begin{lemma}[Projection Lemma]
Let $ \Psi \in \Re^{n \times n}$ be a symmetric matrix and $ P, Q \in \Re^{m \times n}$. Then, following statements are equivalent:%
\begin{enumerate}[label=(\roman*)]
\item There exists $ X \in \Re^{m \times m} $ such that $ \Psi + P^T X^T Q + Q^T X P \preceq 0 $;
\item $ N_P^T \Psi N_P \preceq 0 $ and $ N_Q^T \Psi N_Q \preceq 0$,
\end{enumerate}%
where $ N_P $ and $ N_Q $ are arbitrary bases of the null-space of $ P $ and $ Q $, respectively.

\begin{proof}
See Lemma 3.1 of \citep{gahinet1994linear}.
\end{proof}
\label{lem_projection_lemma}
\end{lemma}

\begin{lemma}[Reciprocal Projection Lemma]
Let $ \Psi \in \Re^{n \times n} $ be a symmetric matrix, $ S \in \Re^{n \times n} $ and $ X \in \Re^{n \times n} $ be an arbitrary positive definite matrix. Then, the following statements are equivalent:%
\begin{enumerate}[label=(\roman*)]
\item $ \Psi + S + S^T \preceq 0 $;
\item There exists $ W \in \Re^{n \times n} $ such that%
\begin{equation}
\begin{bmatrix}
- X & S + W\\
S^T + W^T & \Psi + X - W - W^T 
\end{bmatrix} \preceq 0.
\end{equation}
\end{enumerate}

\begin{proof}
See Lemma 2.2 of \citep{apkarian2001continuous}.
\end{proof}
\label{lem_reciprocal_projection_lemma}
\end{lemma}

\subsection{Property of Kronecker products}

We also present a fundamental property of Kronecker product, necessary for the proposed controller proof.

\begin{lemma}
Let $ X \in \Re^{p \times p} $ and $ Y \in \Re^{q \times r} $. Then,%
\begin{equation}
Z = (X \otimes I_q) (I_p \otimes Y) = (I_p \otimes Y) (X \otimes I_r).
\end{equation}

\begin{proof}
This property follows directly from the fact that each block of $ Z $ is given by $ Z_{i,j} = X_{i,j} Y $, which is commutable since $ X_{i,j} $ is a scalar.
\end{proof}
\label{lem_kronercker_swap}
\end{lemma}

\begin{corollary}
For all $ i \in [1, n] $, let $ X_i \in \Re^{p_i \times p_i} $, $ Y_i \in \Re{q_i \times r_i} $, and%
\begin{equation}
Y = \text{diag}(I_{p_1} \otimes Y_1, I_{p_2} \otimes Y_2, \ldots, I_{p_n} \otimes Y_n).
\end{equation}%
Then, %
\begin{multline}
\text{diag}(X_1 \otimes I_{q_1}, X_2 \otimes I_{q_2}, \ldots, X_n \otimes I_{q_n}) Y = Y \text{diag}(X_1 \otimes I_{r_1}, X_2 \otimes I_{r_2}, \ldots, X_n \otimes I_{r_n}) 
\end{multline}%
follows directly from Lemma \ref{lem_kronercker_swap}.
\label{cor_diag_kronercker_swap}
\end{corollary}

\section{Guaranteed cost control for systems with structured uncertainties}
\label{sec_gcc}

In this section, we present the proposed static output-feedback controller. Existence condition for the particular case of no measured disturbance $ D_y^w = 0 $ is first presented and we specialize it to the state-feedback controller. Then, we present the Optimal Static Output-Feedback Guaranteed Cost Controller developed.

\begin{lemma}[Unstructured Uncertainty Set]
Let $ \xi_k = [x_k^T, w_k^T]^T $. Then, the region where $ w_k = \Delta_k z_k $ holds for some $ \Delta_k \in \mathbb{D}^u $ is%
\begin{equation}
\xi_k \in\mathbb{X}^u = \{ \xi \mid \xi \in \Re^{n_x + n_p}, \xi^T S \xi \le 0, \forall S \in \mathbb{S}^u \},
\end{equation}%
where $ \mathbb{S}^u = \{ \lambda S^u \mid \lambda \ge 0 \} $ is the unstructured uncertainty set, and%
\begin{equation}
S^u = \begin{bmatrix}
\bar{C}_z^T \bar{C}_z & \bar{C}_z^T \bar{D}_z^w\\
\bar{D}_z^{wT} \bar{C}_z & \bar{D}_z^{wT} \bar{D}_z^w - I
\end{bmatrix}.
\end{equation}

\begin{proof}
Direct from the discrete-time Bounded Real Lemma \citep{de1992discrete}.
\end{proof}
\label{lem_disturbance_set_unstructured}
\end{lemma}

The set $ \mathbb{S}^u $, from Lemma \ref{lem_disturbance_set_unstructured}, has been thoroughly used for the analysis and synthesis of robust controllers. It defines the region of the domain of $ \xi_k $ where we wish to ensure system stability and performance. In the case considered in this paper, the structure of the uncertainty $ \Delta_k $ is known and can be used to provide a less conservative set, which still preserves the losslessness property from $ \mathbb{S}^u $.

\begin{lemma}[Structured Uncertainty Set]
Let $ \Lambda_p $, and $ \Lambda_q $ be%
\begin{equation}
\begin{aligned}
\Lambda_p = \text{diag}(\Lambda_1 \otimes I_{n_{p1}}, \Lambda_2 \otimes I_{n_{p2}}, \ldots, \Lambda_s \otimes I_{n_{ps}}),\\
\Lambda_q = \text{diag}(\Lambda_1 \otimes I_{n_{q1}}, \Lambda_2 \otimes I_{n_{q2}}, \ldots, \Lambda_s \otimes I_{n_{qs}}),
\end{aligned}
\end{equation}%
where for all $ 1 \le i \le s $, $ \Lambda_i \in \Re^{n_{ri} \times n_{ri}} $ are positive semi-definite matrices. Then, the region where $ w_k = \Delta_k z_k $ holds for some $ \Delta_k \in \mathbb{D} $ is%
\begin{equation}
\xi_k \in \mathbb{X} = \{ \xi \mid \xi \in \Re^{n_x + n_p}, \xi^T S \xi \le 0, \forall S \in \mathbb{S} \},
\end{equation}%
where $ \mathbb{S} = \{S^s(\Lambda_p, \Lambda_q) \mid \forall i \in [1, s], \forall \Lambda_i \succeq 0 \} $ is the structured uncertainty set, and%
\begin{equation}
S^s(\Lambda_p, \Lambda_q) = \begin{bmatrix}
\bar{C}_z^T \Lambda_q \bar{C}_z & \bar{C}_z^T \Lambda_q \bar{D}_z^w\\
\bar{D}_z^{wT} \Lambda_q \bar{C}_z & \bar{D}_z^{wT} \Lambda_q \bar{D}_z^w - \Lambda_p
\end{bmatrix}.
\end{equation}

\begin{proof}
By left-multiplying the feedback disturbance equality $ w_k = \Delta_k z_k $ by $ \Lambda_p^{\frac{1}{2}} $, we obtain $ \Lambda_p^{\frac{1}{2}} w_k = \Lambda_p^{\frac{1}{2}} \Delta_k z_k $. Which, based on Corollary \ref{cor_diag_kronercker_swap}, is equivalent to%
\begin{equation}
\Lambda_p^{\frac{1}{2}} w_k = \Delta_k \Lambda_q^{\frac{1}{2}} z_k. 
\label{eq_structured_uncertainty_set_proof_1}
\end{equation}%
Consider the 2-norm of \eqref{eq_structured_uncertainty_set_proof_1}. Then, for all $ \Delta_k \in \mathbb{D} $,%
\begin{equation}
|| \Lambda_p^{\frac{1}{2}} w_k ||_2 = ||\Delta_k \Lambda_q^{\frac{1}{2}} z_k||_2 \le ||\Delta_k||_2 ||\Lambda_q^{\frac{1}{2}} z_k||_2 \le ||\Lambda_q^{\frac{1}{2}} z_k||_2,
\label{eq_structured_uncertainty_set_proof_2}
\end{equation}%
which is equivalent to $ w_k^T \Lambda_p w_k - z_k^T \Lambda_q z_k \le 0 $, or $ \xi_k^T S^s(\Lambda_p, \Lambda_q) \xi_k \le 0$. Therefore, \eqref{eq_structured_uncertainty_set_proof_2} holds for all $ \Lambda_p \succ 0 $, $ \Lambda_q \succ 0 $, and $ \Delta_k \in \mathbb{D} $, and the domain of $ \xi_k $ with admissible feedback disturbance is defined by%
\begin{equation}
\xi_k \in \mathbb{X} = \{ \xi \mid \xi \in \Re^{n_x + n_p}, \xi^T S \xi \le 0, \forall S \in \mathbb{S} \}.
\end{equation}
\end{proof}
\label{lem_structured_uncertainty_set}
\end{lemma}

\begin{lemma}
The Structured Uncertainty Set has the following properties:%
\begin{enumerate}[label=(\roman*)]
\item $ \mathbb{S}^u \subseteq \mathbb{S} $;
\item The set $ \mathbb{S} $ is lossless.
\end{enumerate}

\begin{proof}
See Appendix \ref{appendix_lossless_proof}.
\end{proof}
\label{lem_structured_uncertainty_set_properties}
\end{lemma}

We are now ready to investigate the necessary and sufficient conditions for the existence of a GCC that satisfies Definition \ref{def_gcc}.

\begin{theorem}
Consider a static output-feedback controller $ u_k = -K y_k $ regulating the System \eqref{eq_sys_model} subject to structured uncertainties of the form \eqref{eq_uncertainty_structure}. Then, the following statements are equivalent:%
\begin{enumerate}[label=(\roman*)]
\item The controller satisfies Definition \ref{def_gcc};
\item There exists a value function $ V(x_k) = x_k P x_k $, where $ P \succ 0 $, such that for all $ \Delta_k \in \mathbb{D} $%
\begin{equation}
V(x_{k+1}) - V(x_k) \le - x_k^T \left( \bar{C}_c + \bar{D}_c^w \bar{\Delta}_k \bar{C}_z \right)^T (\bullet) x_k;
\label{eq_bellman_gcc}
\end{equation}
\item There exists $ P \succ 0 $ and $ \Lambda_i \succeq 0 $ for $ i \in [1, s] $ such that%
\begin{multline}
\begin{bmatrix}
\bar{A}^T P \bar{A} - P & \bar{A}^T P \bar{B}^w\\ 
\bar{B}^{wT} P \bar{A} & \bar{B}^{wT} P \bar{B}^w
\end{bmatrix}
+
\begin{bmatrix}
\bar{C}_c^T \bar{C}_c & \bar{C}_c^T \bar{D}_c^w\\ 
\bar{D}_c^{wT} \bar{C}_c & \bar{D}_c^{wT} \bar{D}_c^{wT}
\end{bmatrix}
+ \\ +
\begin{bmatrix}
\bar{C}_z^T \Lambda_q \bar{C}_z & \bar{C}_z^T \Lambda_q \bar{D}_z^w\\ 
\bar{D}_z^{wT} \Lambda_q \bar{C}_z & \bar{D}_z^{wT} \Lambda_q \bar{D}_z^{wT} - \Lambda_p
\end{bmatrix}
\preceq 0.
\label{eq_gcc_lmi_1}
\end{multline}
\end{enumerate}

\begin{proof}
We first prove the equivalnece of \textit{(i)} and \textit{(ii)}, then we prove the equivalence between \textit{(ii)} and \textit{(iii)}.

\textit{(i)} $\Rightarrow$ \textit{(ii)}: Let $ c(x_k, u_k) = x_k^T Q x_k + 2 x_k^T N u_k + u_k^T R u_k $ for brevity and consider the value function $ V^*_k(x_k) $. From Bellman's optimality principle, we obtain%
\begin{equation}
V^*_k(x_k) = \underset{\Delta_k \in \mathbb{D}}{\max} \; \underset{u_k}{\min} \; c(x_k, u_k) + V^*_{k+1}(x_{k+1}).
\label{eq_gcc_lmi_1_proof_1}
\end{equation}%
In the infinite horizon case, \eqref{eq_gcc_lmi_1_proof_1} becomes%
\begin{equation}
V^*(x_k) = \underset{\Delta_k \in \mathbb{D}}{\max} \; \underset{u_k}{\min} \; c(x_k, u_k) + V^*(x_{k+1}).
\label{eq_gcc_lmi_1_proof_1_2}
\end{equation}

Assume \textit{(i)} holds, then from Definition \ref{def_gcc} there exists a stabilizing controller $ u_k = -K y_k $ and a sub-optimal value function $ V(x_k) = x_k^T P x_k \ge V^*(x_k) $ such that%
\begin{equation}
V(x_k) = \underset{\Delta_k \in \mathbb{D}}{\max} \; c(x_k, -K y_k) + V(x_{k+1})
\end{equation}
which implies%
\begin{equation}
\forall \Delta_k \in \mathbb{D}: V(x_k) \ge c(x_k, -K y_k) + V(x_{k+1}).
\label{eq_gcc_lmi_1_proof_2}
\end{equation}

Therefore, \eqref{eq_gcc_lmi_1_proof_2} holds for the quadratic form and is equivalent to \textit{(ii)}, since $ c(x_k, -K y_k) = x_k^T \left( \bar{C}_c + \bar{D}_c^w \bar{\Delta}_k C_z \right)^T (\bullet) x_k $.

\textit{(i)} $\Leftarrow$ \textit{(ii)}: Now assume \textit{(i)} does not holds, then there isn't a sub-optimal value function with quadratic form that upper-bounds \eqref{eq_gcc_lmi_1_proof_1_2}. In such case, the infinite horizon Bellman's optimality principle does not have a solution. Therefore, \textit{(ii)} also doesn't hold.

\textit{(ii)} $\Rightarrow$ \textit{(iii)}: Assume \textit{(ii)} holds and consider the system model \eqref{eq_sys_model_cl}, then \eqref{eq_gcc_lmi_1_proof_2} is equivalent to $ \xi_k^T M \xi_k \le 0 $, where
\begin{equation}
M = \begin{bmatrix}
\bar{A}^T P \bar{A} - P & \bar{A}^T P \bar{B}^w\\ 
\bar{B}^{wT} P \bar{A} & \bar{B}^{wT} P \bar{B}^w
\end{bmatrix}
+
\begin{bmatrix}
\bar{C}_c^T \bar{C}_c & \bar{C}_c^T \bar{D}_c^w\\ 
\bar{D}_c^{wT} \bar{C}_c & \bar{D}_c^{wT} \bar{D}_c^{wT}
\end{bmatrix}
\label{eq_gcc_lmi_1_proof_3}
\end{equation}%
for all $ \Delta_k \in \mathbb{D} $. From Lemma \ref{lem_structured_uncertainty_set}, \eqref{eq_gcc_lmi_1_proof_3} holds for all $ \xi \in \mathbb{X} $. Then, based on the generalized S-Procedure, there exists $ S \in \mathbb{S} $ such that $ M - S \preceq 0 $. Or equivalently, there exists $ \Lambda_i \succeq 0 $ for $ i \in [1, s] $ such that \eqref{eq_gcc_lmi_1} holds. Therefore, \textit{(iii)} holds.

\textit{(ii)} $\Leftarrow$ \textit{(iii)}: Conversely, if \textit{(ii)} does not hold, there doesn't exist a $ P \succ 0 $ such that \eqref{eq_gcc_lmi_1_proof_3} holds. Therefore \textit{(iii)} also does not hold.
\end{proof}
\label{the_gcc}
\end{theorem}

\subsection{Convex condition for systems without disturbance feed-through}

If $ D_y^w = 0 $ is assumed, we obtain $ \bar{B}^w = B^w $, $ \bar{D}_z^w = D_z^w $, and $ \bar{D}_c^w = \mathbf{0_{n_c \times n_p}} $. In such a case, a simplified existence condition can be obtained.

\begin{lemma}
Assume $ D_y^w = 0 $. Then, a static output-feedback controller $ u_k = - K y_k $ is said to be of guaranteed cost, according to Definition \ref{def_gcc}, if and only if there exists $ X \succ 0 $ and $ \Upsilon_i \succeq 0 $ for $ i \in [1, s] $ such that%
\begin{equation}
\begin{bmatrix}
- \Upsilon_q & 0 & 0 & C_z X - D_z^u Y C_y & D_z^w \Upsilon_p\\ 
\star & - I_{n_c} & 0 & C_c X - D_c^u Y C_y & 0\\ 
\star & \star & -X & A X - B^u Y C_z & B^w \Upsilon_p\\ 
\star & \star & \star & -X & 0\\ 
\star & \star & \star & \star & - \Upsilon_p
\end{bmatrix} \preceq 0,
\label{eq_gcc_lmi_2}
\end{equation}%
where $ \bar{X} C_y = C_y X $ and $ Y = K \bar{X} $.

\begin{proof}
We first apply the Schur complement in \eqref{eq_gcc_lmi_1} to the terms dependent on $ P $, to the cost terms, and to the terms dependent on $ \Lambda_q $, which results%
\begin{equation}
\begin{bmatrix}
- \Lambda_q^{-1} & 0 & 0 & \bar{C}_z & D_z^w\\ 
\star & - I_{n_c} & 0 & \bar{C}_c & 0\\ 
\star & \star & -P^{-1} & \bar{A} & B^w\\ 
\star & \star & \star & -P & 0\\ 
\star & \star & \star & \star & - \Lambda_p
\end{bmatrix} \preceq 0.
\label{eq_gcc_lmi_2_proof_1}
\end{equation}%
Then, we apply the congruence transformation in \eqref{eq_gcc_lmi_2_proof_1} with $ T = \text{diag}(I_{n_q}, I_{n_c}, I_{n_x}, P^{-1}, \Lambda_p^{-1}) $, and perform the substitution $ X = P^{-1} $, $ \Upsilon_p = \Lambda_p^{-1} $, and $ \Upsilon_q = \Lambda_q^{-1} $, obtaining%
\begin{equation}
\begin{bmatrix}
- \Upsilon_q & 0 & 0 & \bar{C}_z X & D_z^w \Upsilon_p\\ 
\star & - I_{n_c} & 0 & \bar{C}_c X & 0\\ 
\star & \star & -X & \bar{A} X & B^w \Upsilon_p\\ 
\star & \star & \star & -X & 0\\ 
\star & \star & \star & \star & - \Upsilon_p
\end{bmatrix} \preceq 0.
\label{eq_gcc_lmi_2_proof_2}
\end{equation}%
Based on the substitutions $ \bar{X} C_y = C_y X $, $ Y = K \bar{X} $, $ \bar{A} X = A X - B^u Y C_z$, $ \bar{C}_c X = C_c X - D_c^u Y C_y $, and $ \bar{C}_z X = C_z X - D_z^u Y C_y $, \eqref{eq_gcc_lmi_2} and \eqref{eq_gcc_lmi_2_proof_2} are equivalent.
\end{proof}
\label{lem_gcc_equiv_1}
\end{lemma}

\begin{remark}
The control gain can be recovered from Lemma \ref{lem_gcc_equiv_1} by $ K = Y (C_y X C_y^\dagger)^{-1} $, where $ C_y^\dagger $ is the pseudo-inverse of $ C_y $. 
\end{remark}

\begin{remark}
Given $ C_y = I_{n_x} $, we obtain $ X = \bar{X} $, $ K = Y X^{-1} $, and Lemma \ref{lem_gcc_equiv_1} reduces to the state-feedback controller synthesis condition.
\end{remark}

Lemma \ref{lem_gcc_equiv_1} provides an LMI condition for the existence of the GCC, based on the multipliers $ \Lambda_i $ for $ i \in [1, s] $, the inverse cost function matrix $ X $, and $ Y $. From the properties of LMIs we can conclude that \eqref{eq_gcc_lmi_2} is convex.

\subsection{Convex condition for systems with disturbance feed-through}

If we consider the general case, where $ D_y^w \neq 0 $, the controller gain $ K $ would also be present in the fifth column of \eqref{eq_gcc_lmi_2}. In such a case, the substitution which yields $ Y $, used in Lemma \ref{lem_gcc_equiv_1}, would not make the condition linear since $ K $ would still exist in the fifth row and column.

The following Theorem is the main result of this paper. It demonstrates that a similar substitution can be successfully used to convexify a ``dilated" version of condition \eqref{eq_gcc_lmi_2}, based on the Reciprocal Projection Lemma \ref{lem_reciprocal_projection_lemma}.

\begin{theorem}
A static output-feedback controller $ u_k = - K y_k $ is said to be of guaranteed cost, according to Definition \ref{def_gcc}, if there exists $ X \succ 0 $, $ \Lambda_i \succeq 0 $ for $ i \in [1, s] $, $ V_{2,2} \in \Re^{n_c \times n_c} $, $ V_{3,3} \in \Re^{n_x \times n_x} $, $ V_{4,4} \in \Re^{n_x \times n_x} $, $ V_{4,5} \in \Re^{n_x \times n_q} $, $ V_{5,4} \in \Re^{n_q \times n_x} $, $ V_{5,5} \in \Re^{n_q \times n_q} $, and $ Y \in \Re^{n_u \times n_y} $ such that%
\begin{equation}
\begin{bmatrix}
- M & 0 & V\\
\star & -M & \bar{N} + M\\
\star & \star & - V - V^T
\end{bmatrix} \preceq 0
\label{eq_gcc_lmi_4}
\end{equation}%
where, $ M = \text{diag}(\Upsilon_q, I_{n_c}, X, X, \Upsilon_q) $,%
\begin{equation}
V = \begin{bmatrix}
V_{1,1} & V_{1,2} & V_{1,3} & V_{1,4} & V_{1,5}\\
V_{2,1} & V_{2,2} & V_{2,3} & V_{2,4} & V_{2,5}\\
V_{3,1} & V_{3,2} & V_{3,3} & V_{3,4} & V_{3,5}\\
0 & 0 & 0 & V_{4,4} & V_{4,5}\\
0 & 0 & 0 & V_{5,4} & V_{5,5}\\
\end{bmatrix},
\label{eq_gcc_lmi_4_v}
\end{equation}%
\begin{equation}
\bar{N} =
\frac{1}{2}
\begin{bmatrix}
- V_{1,1} & - V_{1,2} & - V_{1,3} & \Phi_z & \Sigma_z\\
- V_{2,1} & - V_{2,2} & - V_{2,3} & \Phi_c & \Sigma_c\\
- V_{3,1} & - V_{3,2} & - V_{3,3} & \Phi_x & \Sigma_x\\
0 & 0 & 0 & - V_{4,4} & - V_{4,5}\\
0 & 0 & 0 & - V_{5,4} & - V_{5,5}
\end{bmatrix},
\end{equation}%
\begin{equation}
\begin{bmatrix}
\Phi_z & \Sigma_z\\
\Phi_c & \Sigma_c\\
\Phi_x & \Sigma_x
\end{bmatrix}
=
2 \begin{bmatrix}
C_z & D_z^w \\
C_c & 0 \\
A & B^u \\
\end{bmatrix}
\begin{bmatrix}
V_{4,4} & V_{4,5}\\
V_{5,4} & V_{5,5}
\end{bmatrix}
+ 
\begin{bmatrix}
V_{1,4} & V_{1,5} \\
V_{2,4} & V_{2,5} \\
V_{3,4} & V_{3,5} \\
\end{bmatrix}
-
2 \begin{bmatrix}
D_z^u \\ D_c^u \\ B^u
\end{bmatrix} Y \begin{bmatrix}
C_y & D_y^w
\end{bmatrix},
\end{equation}%
\begin{equation}
\bar{V} \begin{bmatrix}
C_y & D_y^w
\end{bmatrix}
=
\begin{bmatrix}
C_y & D_y^w
\end{bmatrix}
\begin{bmatrix}
V_{4,4} & V_{4,5} \\
V_{5,4} & V_{5,5}
\end{bmatrix},
\end{equation}%
$ Y = K \bar{V} $, $ X = P^{-1}$, $ \Upsilon_p = \Lambda_p^{-1} $, and $ \Upsilon_q = \Lambda_q^{-1} $.

\begin{proof}
We apply again Schur complement in \eqref{eq_gcc_lmi_1} to obtain%
\begin{equation}
\begin{bmatrix}
- \Lambda_q^{-1} & 0 & 0 & \bar{C}_z & \bar{D}_z^w\\ 
\star & - I_{n_c} & 0 & \bar{C}_c & \bar{D}_c^w\\ 
\star & \star & -P^{-1} & \bar{A} & \bar{B}^w\\ 
\star & \star & \star & -P & 0\\ 
\star & \star & \star & \star & - \Lambda_p
\end{bmatrix} \preceq 0.
\label{eq_gcc_lmi_4_proof_1}
\end{equation}%
Then, we apply the congruence transformation in \eqref{eq_gcc_lmi_4_proof_1} with $ T = \text{diag}(\lambda_q, I_{n_c}, P, I_{n_x}, I_{n_p}) $, and obtain%
\begin{equation}
E = \begin{bmatrix}
- \Lambda_q & 0 & 0 & \Lambda_q \bar{C}_z & \Lambda_q \bar{D}_z^w\\ 
\star & - I_{n_c} & 0 & \bar{C}_c & \bar{D}_c^w\\ 
\star & \star & -P & P \bar{A} & P \bar{B}^w\\ 
\star & \star & \star & -P & 0\\ 
\star & \star & \star & \star & - \Lambda_p
\end{bmatrix} \preceq 0.
\label{eq_gcc_lmi_4_proof_2}
\end{equation}

We now define $ S $ as%
\begin{equation}
S = \begin{bmatrix}
- \frac{1}{2} \Lambda_q & 0 & 0 & \Lambda_q \bar{C}_z & \Lambda_q \bar{D}_z^w\\ 
0 & - \frac{1}{2} I_{n_c} & 0 & \bar{C}_c & \bar{D}_c^w\\ 
0 & 0 & - \frac{1}{2} P & P \bar{A} & P \bar{B}^w\\ 
0 & 0 & 0 & - \frac{1}{2} P & 0\\ 
0 & 0 & 0 & 0 & - \frac{1}{2} \Lambda_p
\end{bmatrix},
\label{eq_gcc_lmi_4_proof_3}
\end{equation}%
such that $ S + S^T = E $. Then, from the Reciprocal Projection Lemma (Lemma \ref{lem_reciprocal_projection_lemma}) we obtain that for any given $ Y \succ 0 $ there exists $ W $, of appropriate dimensions, such that%
\begin{equation}
\begin{bmatrix}
- Y & S + W\\
\star & Y - W - W^T 
\end{bmatrix} \preceq 0.
\label{eq_gcc_lmi_4_proof_4}
\end{equation}

Subsequently, we apply a congruence transformation in \eqref{eq_gcc_lmi_4_proof_4} with $ T = \text{diag}(M, I_{n_q+n_c+n_x+n_x+n_p}) $ and substitute $ Y = M^{-1} $, which yields%
\begin{equation}
\begin{bmatrix}
- M & MS + MW\\
\star & M^{-1} - W - W^T 
\end{bmatrix} \preceq 0.
\label{eq_gcc_lmi_4_proof_5}
\end{equation}%
Let $ V = W^{-1} $ and assume it takes the form \eqref{eq_gcc_lmi_4_v}. Then, we apply a congruence transformation in \eqref{eq_gcc_lmi_4_proof_4} with $ T = \text{diag}(I_{n_q+n_c+n_x+n_x+n_p}, V) $, resulting%
\begin{equation}
\begin{bmatrix}
- M & MSV + M\\
\star & V M^{-1} V - V - V^T 
\end{bmatrix} \preceq 0
\label{eq_gcc_lmi_4_proof_6}
\end{equation}%
which, with a Schur complement of the terms related to $ M^{-1} $ and $ MSV = \bar{N} $, is equivalent to%
\begin{equation}
\begin{bmatrix}
- M & 0 & V\\
\star & - M & \bar{N} + M\\
\star & \star & - V - V^T 
\end{bmatrix} \preceq 0
\label{eq_gcc_lmi_4_proof_7}
\end{equation}%
which, in turn, is identical to \eqref{eq_gcc_lmi_4}.
\end{proof}
\label{the_gcc_equiv_3}
\end{theorem}

Theorem \ref{the_gcc_equiv_3} provides a convex condition for the existence of the GCC in the general static output-feedback case. It enables the synthesis of Optimal GCCs through SDP optimization problems.

\subsection{Optimal static Output-feedback guaranteed cost control synthesis}

In this sub-section, we present the resulting optimization problem for the optimal static output feedback GCC synthesis. Although the formulated GCC conditions are state-independent, the cost function \eqref{eq_cost_function} is still dependent. Therefore, we must represent the system cost in a state-independent manner. Two methods are often employed to address this issue: minimizing the worst case cost for all states in the unit norm ball \citep{petersen1998optimal}, or the expected cost value for a zero mean and unit covariance initial state distribution \citep{xie1993control}. We have chosen the latter, since the worst case cost only minimizes the largest eigenvalue of $ P $, while the stochastic interpretation minimizes $ \tr(P) $. Another advantage of such an approach is its equivalence to the LQR when the system reduces to the state-feedback case without uncertainties.

\begin{assumption}
The initial state $ x_0 $ is a zero mean random variable with unit covariance.
\label{ass_state_covariance}
\end{assumption}

From Assumption \ref{ass_state_covariance}, $ E(J^*(x_0)) = E(x_0^T P x_0) = E(\tr(x_0^T P x_0)) = \tr(P E(x_0 x_0^T)) = \tr(P) $.

\begin{theorem}
A static output feedback controller $ u_k = - K y_k $ is said to be an optimal guaranteed cost, according to Definition \ref{def_gcc_optimal}, if and only if%
\begin{equation}
\begin{aligned}
K = \argmin \;\; & \tr(Z) \\
s.t. \;\; & \begin{bmatrix}
- M & 0 & V\\
\star & -M & \bar{N} + M\\
\star & \star & - V - V^T
\end{bmatrix} \preceq 0\\
& \begin{bmatrix}
-Z & I\\
\star & -X
\end{bmatrix} \preceq 0
\end{aligned}
\label{eq_optimal_gcc}
\end{equation}%
where $ M $, $\bar{N} $, $ V $ and $ X $ are defined in Theorem \ref{the_gcc_equiv_3}, and $ Z \in \Re^{n_x \times n_x} $ is the cost matrix $ P $ relaxation.

\begin{proof}
From Definition \ref{def_gcc_optimal}, a controller is of optimal guaranteed cost if it satisfies Definition \ref{def_gcc} and $ J^*(x_0) $ is minimal. Since the distribution of $ x_0 $ is assumed to be known (Assumption \ref{ass_state_covariance}), consider $ \widehat{J} = E(J(x_0)) $. Then, the optimal GCC synthesis is equivalent to%
\begin{equation}
\begin{aligned}
\widehat{J}^* = \min \;\; & \tr(P) \\
s.t. \;\; & \begin{bmatrix}
- M & 0 & V\\
\star & -M & \bar{N} + M\\
\star & \star & - V - V^T
\end{bmatrix} \preceq 0
\end{aligned}
\label{eq_optimal_gcc_proof_1}
\end{equation}%
which minimizes the expected cost. However, such optimization is not linear on $ P $ since the LMI depends on its inverse $ X = P^{-1} $. Therefore, we introduce the relaxation variable $ Z \in \Re^{n_x \times n_x} $, such that $ Z \succeq P $. From Schur complement, we obtain%
\begin{equation}
\begin{bmatrix}
-Z & I\\
\star & -X
\end{bmatrix} 
\preceq 0.
\end{equation}%
Notice that minimizing Z is equivalent to minimizing $ P $, since its constraint is only dependent on $ X $. Therefore, \eqref{eq_optimal_gcc_proof_1} is equivalent to \eqref{eq_optimal_gcc}.
\end{proof}
\label{the_optimal_gcc}
\end{theorem}

\begin{remark}
The optimal static output feedback for the case $ D_y^w = 0 $, can be defined similarly to Theorem \ref{the_gcc_equiv_3}, with the substitution of condition \eqref{eq_gcc_lmi_4} by \eqref{eq_gcc_lmi_2}.
\label{rem_other_gccs_sdp}
\end{remark}

\section{Numerical examples}
\label{sec_example}

This section provides two numerical examples of GCC for systems subject to structured uncertainties. The first compares the proposed synthesis method with GCC for systems subject to unstructured uncertainties \cite{xie1993control}, and the Linear Quadratic Regulator \cite{kalman1960contributions}. The second presents the synthesis and simulation for a system that previously proposed methods cannot synthesize a controller. We have used the YALMIP Toolbox \cite{lofberg2004yalmip} and the Mosek solver \cite{mosek} for the modeling of the problem\footnotemark[1].

\footnotetext[1]{The source code for both numerical examples are available at: \url{https://gitlab.com/cmasseraf/gcc_structured_uncertainty}}

\subsection{Example 1}

Consider the uncertain linear system (\eqref{eq_sys_model} and \eqref{eq_sys_uncertainty_matrix}), from \cite{terra2014optimal}, defined by%
\begin{equation}
\begin{matrix}
A = \begin{bmatrix}
1.1 & 0 & 0\\ 
0 & 0 & 1.2\\ 
-1 & 1 & 0
\end{bmatrix},
\;
B^u = \begin{bmatrix}
0 & 1\\ 
1 & 1\\ 
-1 & 0
\end{bmatrix},
\;
B^w = \begin{bmatrix}
0.7 & 0.3\\ 
0.5 & -0.4\\ 
-1 & 0
\end{bmatrix},
\\
C_z = \begin{bmatrix}
0.41 & 0.43 & -0.5\\ 
0 & -0.32 & 0.44
\end{bmatrix},
\;
D_z^u = \begin{bmatrix}
0.4 & -0.4\\ 
0 & 0
\end{bmatrix},
\;
D_z^w = 0,
\end{matrix}
\label{eq_model_example}
\end{equation}%
measurement matrices%
\begin{equation}
C_y = I_{n_x},\;D_y^w = 0.
\end{equation}%
and the disturbance%
\begin{equation}
\Delta_k = \begin{bmatrix}
\delta_{1, k} & 0\\
0 & \delta_{2, k}
\end{bmatrix}
\end{equation}%
where $ \delta_{1,k}, \delta_{1,k} \in [-1, 1] $. In this example, we synthesize and compare controllers for five different methods.

\begin{enumerate}
\item LQR \cite{kalman1960contributions}: The LQR controller is generated based on the nominal plant with no modeled disturbances.

\item GCC (unstructured uncertainty) \cite{xie1993control}: Based on Xie's approach.

\item Theorem \ref{the_gcc_equiv_3}-based GCC (unstructured): Controller based on Theorem \ref{the_gcc_equiv_3}, generated with the SDP optimization of Theorem \ref{the_optimal_gcc}, while assuming no structure to $ \Delta_k $. Therefore, $ \Lambda_p = \lambda I_p $ and $ \Lambda_q = \lambda I_q $ where $ \lambda > 0 $.

\item Lemma \ref{lem_gcc_equiv_1}-based GCC (structured): Controller based on Theorem \ref{lem_gcc_equiv_1}, generated with the SDP optimization of Theorem \ref{the_optimal_gcc} without the ``Dilated" LMI (see Remark \ref{rem_other_gccs_sdp}), while exploiting the diagonal uncertainty structure of $ \Delta_k $. Therefore, $ \Lambda_p = \text{diag}(\lambda_1, \lambda_2) $ and $ \Lambda_q = \text{diag}(\lambda_1, \lambda_2) $ where $ \lambda_i > 0 $ for all $ i \in [1,2] $.

\item Theorem \ref{the_gcc_equiv_3}-based GCC (structured): Controller based on Theorem \ref{the_gcc_equiv_3}, generated with the SDP optimization of Theorem \ref{the_optimal_gcc}, while exploiting the diagonal uncertainty structure.
\end{enumerate}

The cost function matrices selected for all methods were $ Q = I_3 $, and $ R = I_2 $. Different methods are compared based on the synthesis cost, $\tr(P)$, and the effective cost, $ E(J(x_0)) $. We calculated the effective cost as the mean value of the cost function for $5000$ simulations performed with a horizon of $200$ time-steps. Then, we present a simulation example for three controllers synthesized from these methods and qualitatively compare the obtained results.

\begin{figure*}[!t]
\centering
\begin{subfigure}{0.45\textwidth}
 \centering
 \includegraphics[width=\linewidth]{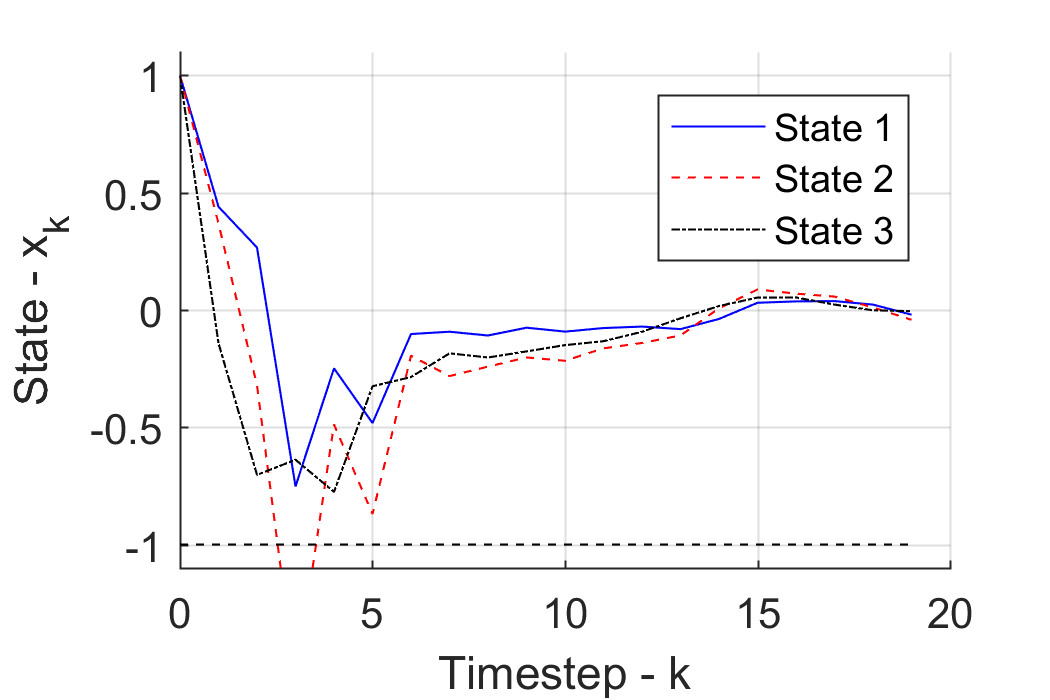}
 \caption{}
 \label{fig_ex1_a}
\end{subfigure}
\quad
\begin{subfigure}{0.45\textwidth}
 \centering
 \includegraphics[width=\linewidth]{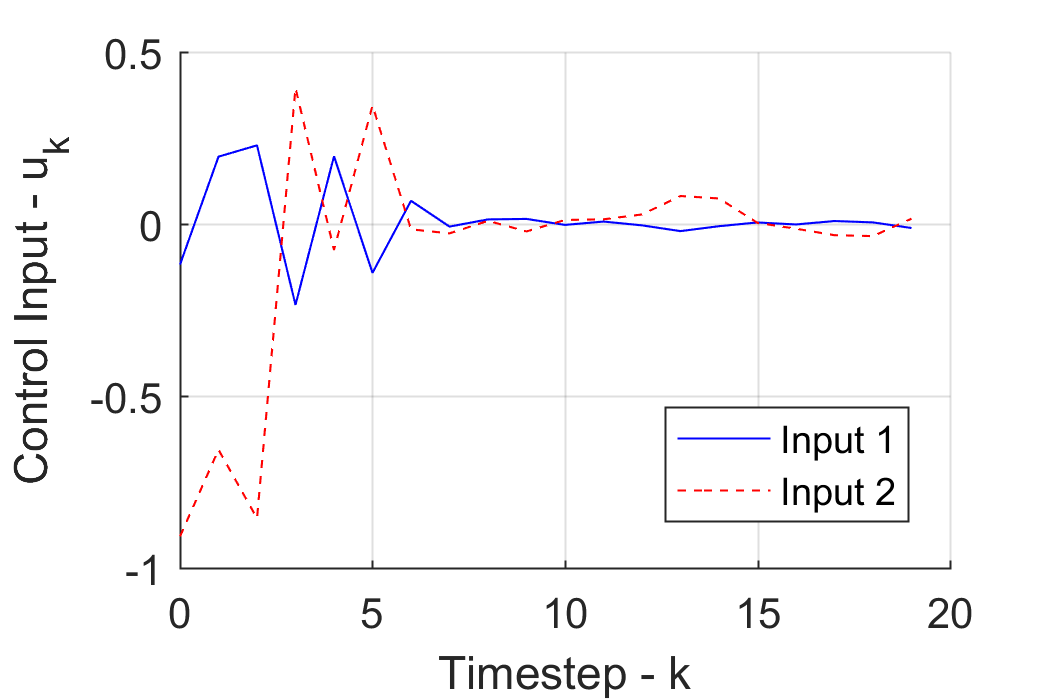}
 \caption{}
 \label{fig_ex1_b}
\end{subfigure}
\begin{subfigure}{0.45\textwidth}
 \centering
 \includegraphics[width=\linewidth]{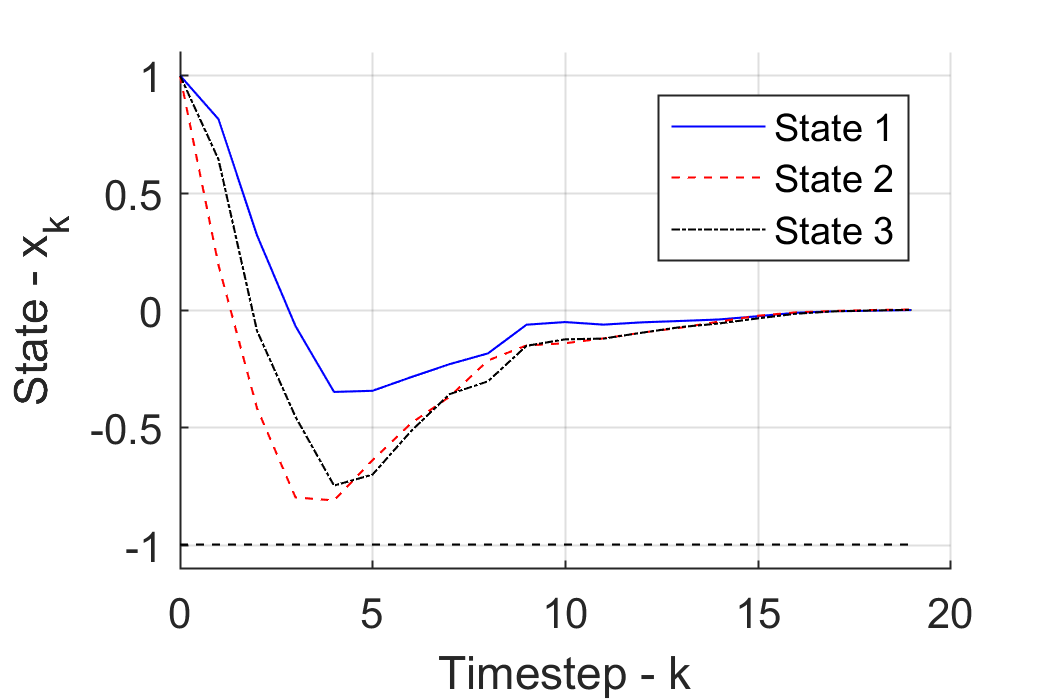}
 \caption{}
 \label{fig_ex1_c}
\end{subfigure}
\quad
\begin{subfigure}{0.45\textwidth}
 \centering
 \includegraphics[width=\linewidth]{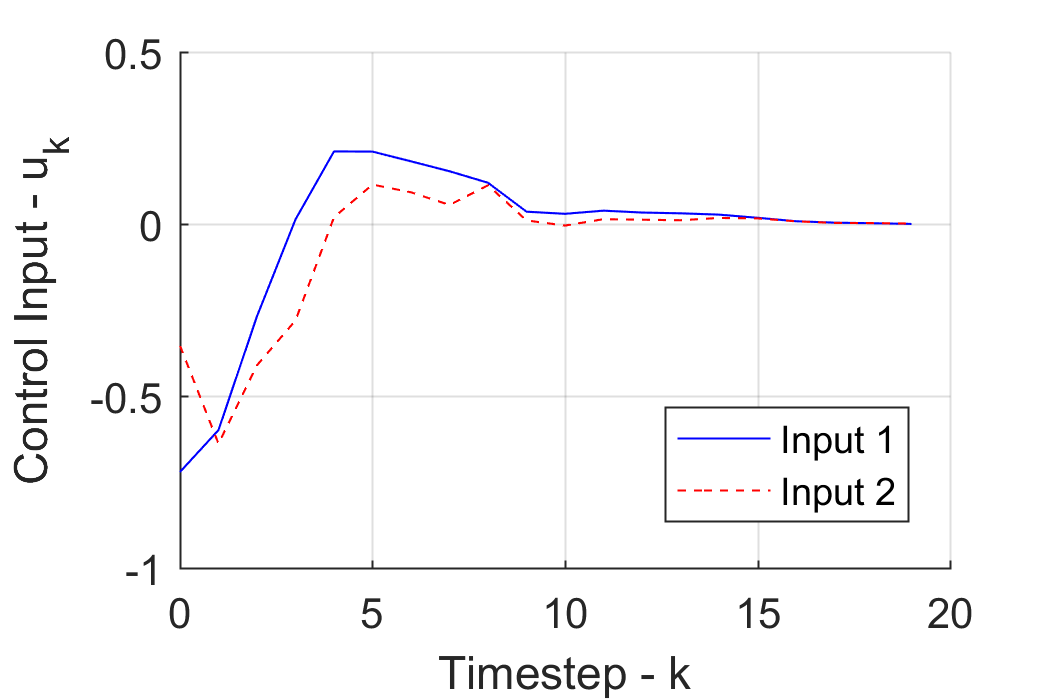}
 \caption{}
 \label{fig_ex1_d}
\end{subfigure}
\begin{subfigure}{0.45\textwidth}
 \centering
 \includegraphics[width=\linewidth]{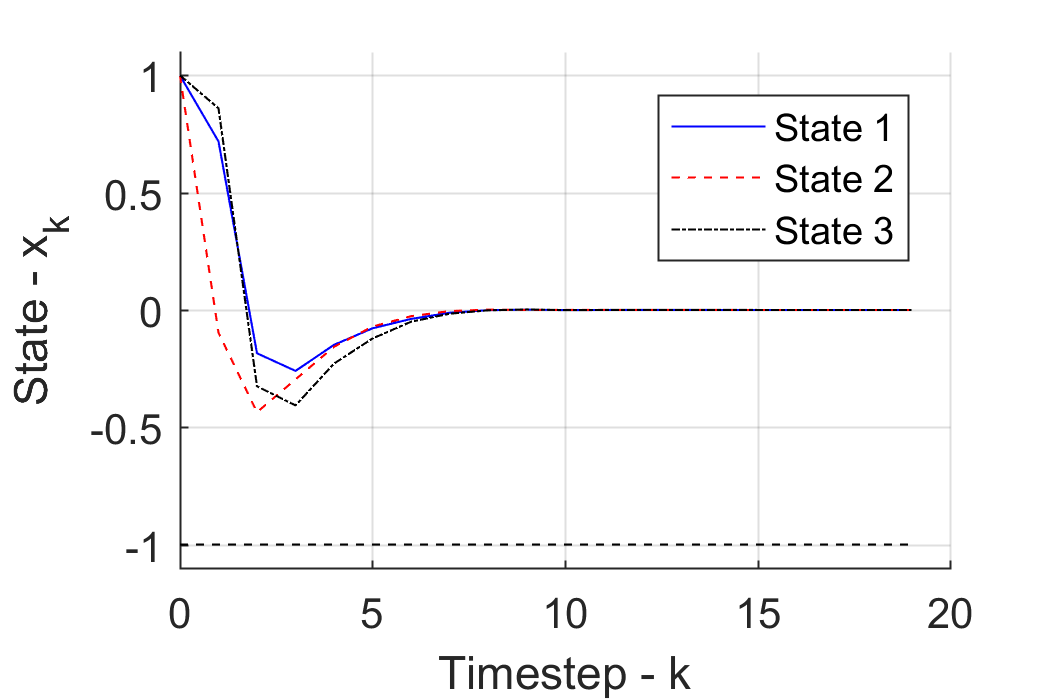}
 \caption{}
 \label{fig_ex1_e}
\end{subfigure}
\quad
\begin{subfigure}{0.45\textwidth}
 \centering
 \includegraphics[width=\linewidth]{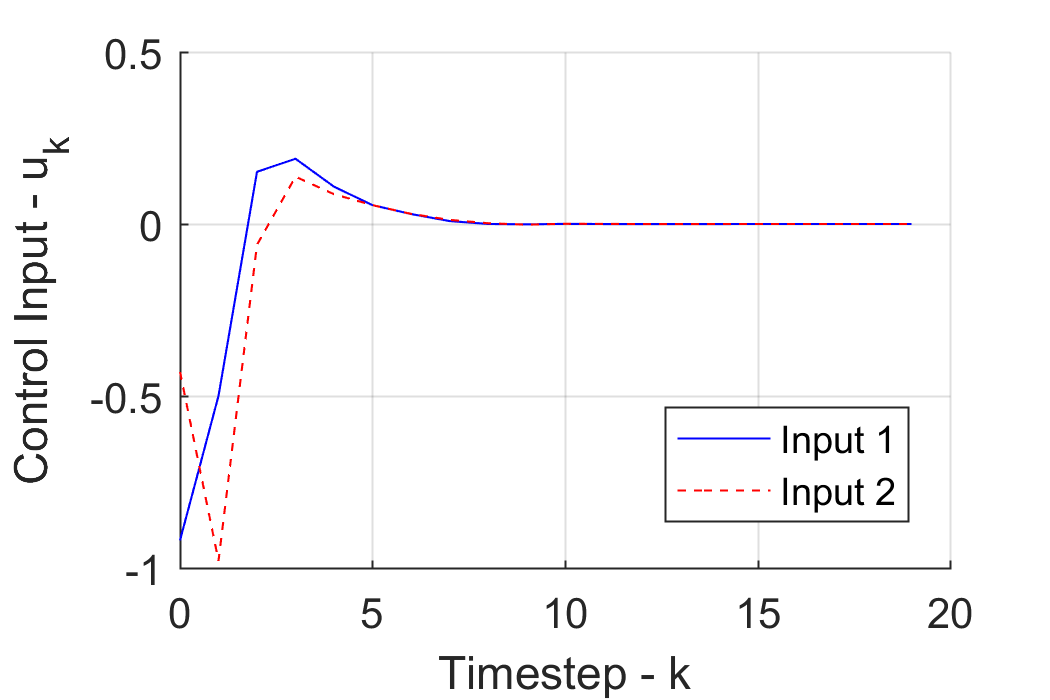}
 \caption{}
 \label{fig_ex1_f}
\end{subfigure}
\caption{Simulation results for the synthesized results of Example 1. (a) Closed-loop state $ x_k $ subject to uncertainties for (1); (b) Control inputs $ u_k $ for (1); (c) Closed-loop state $ x_k $ subject to uncertainties for (2) and (3); (d); Control inputs $ u_k $ for (2) and (3); (e) Closed-loop state $ x_k $ subject to uncertainties for (4) and (5); (f) Control inputs $ u_k $ for (4) and (5).}
\label{fig_example_1_simulation}
\end{figure*}

Table \ref{tbl_example_1} presents the results of Example 1. It is possible to see that for this particular system, the effective cost of the LQR (1) exceeds its synthesis cost by one order of magnitude demonstrating its expected lack of robustness margins, although it is still able to stabilize the system. Both unstructured uncertainty GCCs (2 and 3) yield the same control gain and cost since $ D_y^w = 0 $. Their synthesis cost is greater than the effective cost, which illustrates the upper-bound nature of its robustness margins. However, they yield extremely conservative bounds on cost (of an order of magnitude) due to not exploiting the uncertainty structure. Finally, both structured uncertainty GCCs (4 and 5) also yield the same control gain to each other since $ D_y^w = 0 $, but decrease the conservativeness of the synthesis by a factor of seven while also reducing the effective cost from $46.56$ to $44.16$. 

\begin{table}[ht]
\centering
\caption{Synthesis and simulation results for the comparative study}
\label{tbl_example_1}
\begin{tabular}{l|c|c|}
\cline{2-3}
 & Synthesis cost & Effective cost \\ \hline
\multicolumn{1}{|l|}{1) LQR} & 22.15 & 490.13 \\ \hline
\multicolumn{1}{|l|}{2) GCC proposed by \cite{xie1996output}} & 581.79 & 46.56 \\ \hline
\multicolumn{1}{|l|}{3) Theorem \ref{the_gcc_equiv_3} GCC (unstructured)} & 581.79 & 46.56 \\ \hline
\multicolumn{1}{|l|}{4) Lemma \ref{lem_gcc_equiv_1} GCC (structured)} & 97.63 & 44.16 \\ \hline
\multicolumn{1}{|l|}{5) Theorem \ref{the_gcc_equiv_3} GCC (structured)} & 97.63 & 44.16 \\ \hline
\end{tabular}
\end{table}

Figure 1 presents results acquired by one of the 5000 simulations performed with the set of controllers considered in this example. Identical disturbances and initial states, $ x_0 = [1, 1, 1]^T $, were used for all simulation. Figures \ref{fig_ex1_a} and \ref{fig_ex1_b} show the results for the standard LQR controller with effective cost $ 12.6707 $. We can observe a significant impact of the disturbances on the overall system behavior, due to the lack of robustness guarantees on its synthesis. Figures \ref{fig_ex1_c} and \ref{fig_ex1_d} show the results for the GCC considering unstructured disturbances with effective cost $ 10.7861 $; we can observe that it rejects the disturbances significantly better than the LQR controller; however, it yields an $80\%$ overshoot before reaching steady-state. Finally, Figures \ref{fig_ex1_e} and \ref{fig_ex1_f} show the results for the GCC considering unstructured disturbances with effective cost of $ 7.3934 $. This controller was able to yield a faster settling time, smaller overshoot, and lower effective cost when compared to both other synthesized controllers.

The average cost reduction, combined by the performance improvement observed in Figure 1, demonstrates the advantages of incorporating the uncertainty structure information into the controller synthesis procedure.

\subsection{Example 2}

In this example we present the synthesis and simulation for a system that was not supported by previous methods, where $ D_y^w \neq 0 $. Consider the system from \eqref{eq_model_example} with measurement matrices%
\begin{equation}
C_y = \begin{bmatrix}
1 & 0 & 0\\
0 & 1 & 0\\
0 & 0 & 1\\
0 & 0 & 0
\end{bmatrix},
\;
D_y^w = \begin{bmatrix}
0 & 0\\
0 & 0\\
0 & 0\\
0 & 1\\
\end{bmatrix}.
\end{equation}

In such a case, the LQR (1), Xie's GCC (2) and Lemma \ref{lem_gcc_equiv_1}-based GCC (4) cannot be used as they assume $ D_y^w = 0 $. However, the unstructured uncertainty GCC based on Theorem \ref{the_gcc_equiv_3} (3) also cannot be synthesized, since Assumption \ref{ass_no_feedthrough} ($D_y^w \Delta_k D_z^u = 0$) does not hold for the generic uncertainty case. Therefore, the proposed GCC controller based on Theorem \ref{the_gcc_equiv_3} for structured uncertainties is the only capable of synthesizing a controller for this system.

\begin{figure*}[!t]
\centering
\begin{subfigure}{\textwidth}
 \centering
 \includegraphics[width=\linewidth]{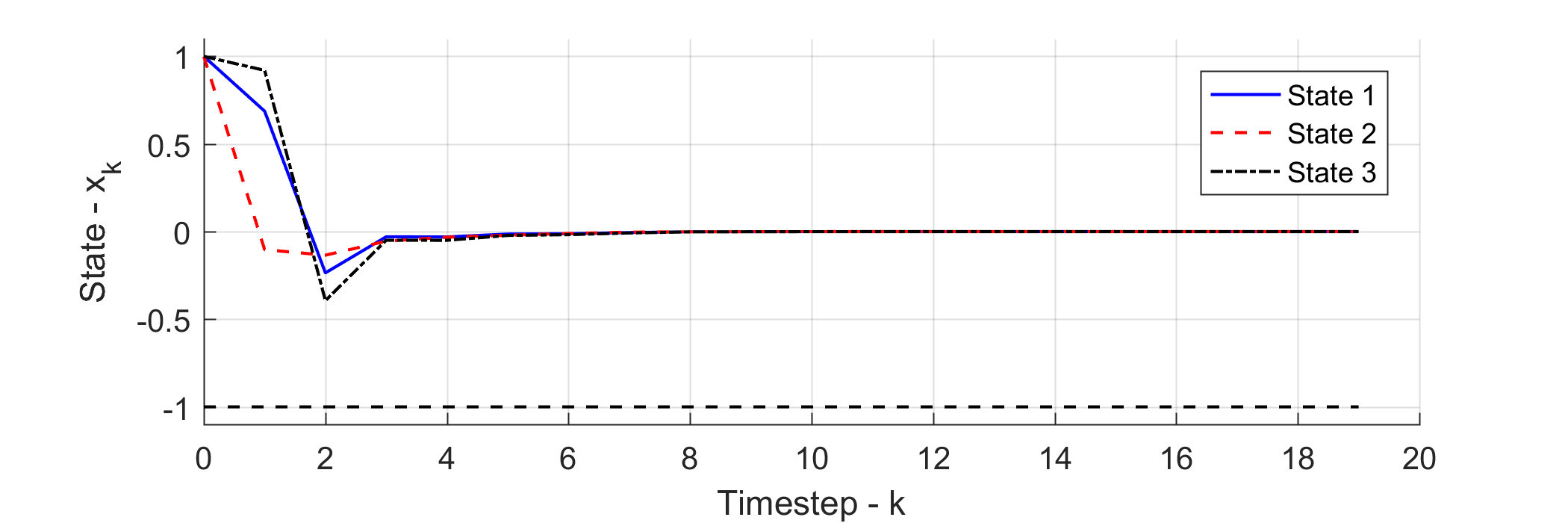}
 \caption{}
 \label{fig_result_b}
\end{subfigure}
\quad
\begin{subfigure}{0.45\textwidth}
 \centering
 \includegraphics[width=\linewidth]{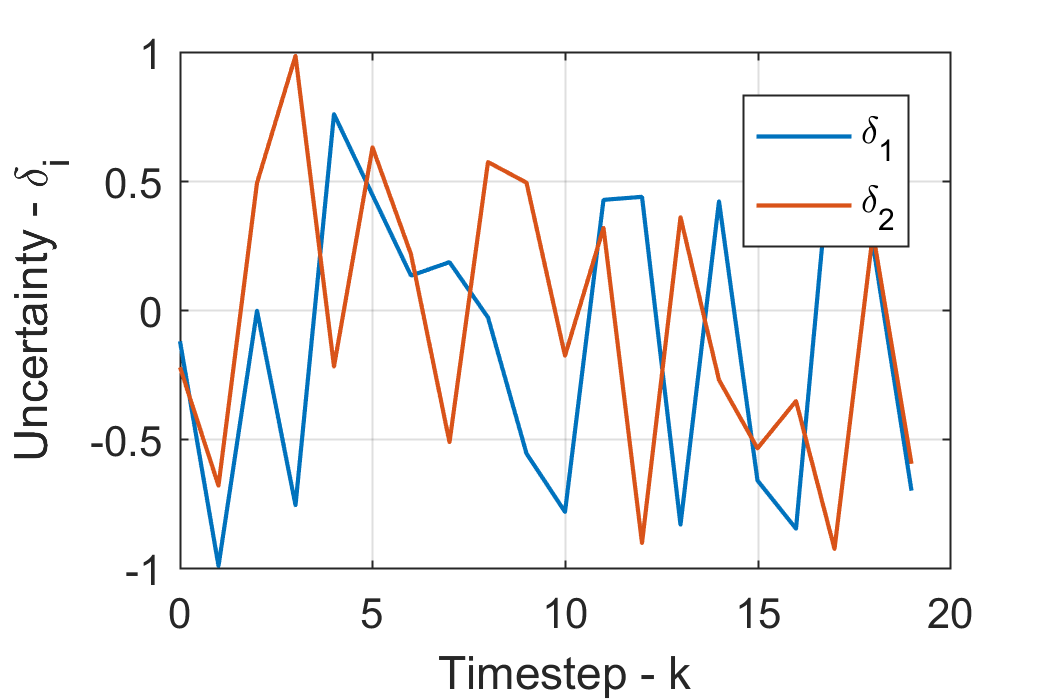}
 \caption{}
 \label{fig_result_a}
\end{subfigure}
\quad
\begin{subfigure}{0.45\textwidth}
 \centering
 \includegraphics[width=\linewidth]{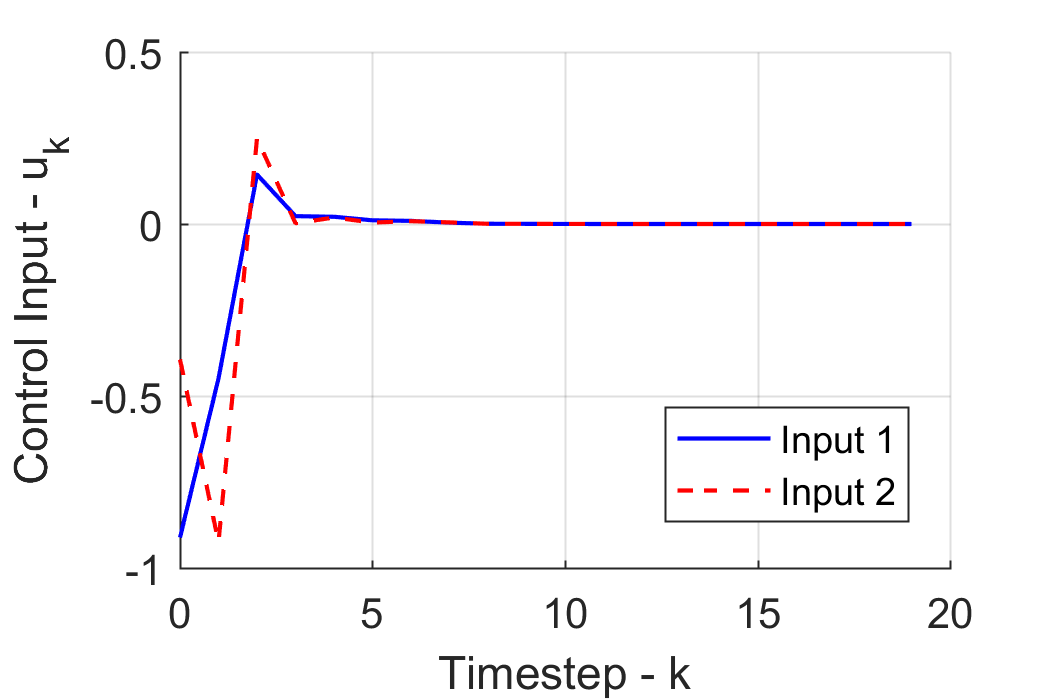}
 \caption{}
 \label{fig_result_c}
\end{subfigure}
\caption{Simulation results for GCC of Example 2. (a) Closed-loop state $ x_k $ subject to uncertainties; (b) Simulated uncertainty $ \delta_1 $ and $ \delta_2 $; (c) Control inputs $ u_k $.}
\label{fig_example_2}
\end{figure*}

The structured uncertainty GCC based on Theorem \ref{the_gcc_equiv_3} synthesis resulted in a gain matrix%
\begin{equation}
K = \begin{bmatrix}
1.1431 & 0.1282 & -0.3585 & 0.0947 \\
0.6881 & -0.7581 & 0.4561 & 0.0596
\end{bmatrix},
\end{equation}%
and cost matrix %
\begin{equation}
P = \begin{bmatrix}
61.9182 & 3.2483 & -41.2619 \\
3.2483 & 9.8246 & -7.1655 \\
-41.2619 & -7.1655 & 34.2462
\end{bmatrix}.
\end{equation}%
It provides a synthesis cost of $94.15$ and an effective cost of $45.62$. Figure \ref{fig_example_2} shows the simulation result for initial state $ x_0 = [1, 1, 1]^T $, where we can observe that the control law is able to correctly stabilize the system subject to varying disturbances.

\section{Conclusion}
\label{sec_conclusion}

In this paper, we have proposed a state-feedback and static output-feedback generalized Guaranteed Cost Control method for discrete-time linear systems subject to norm-bounded structured parametric uncertainties. Such controllers provide upper bounds to the quadratic cost functional and guarantee robust stability and performance to a larger class of uncertain linear systems if compared with previous methods. There are two main contributions for the class of systems considered:
It relaxes the assumption from no uncertainty feed-through ($ D_y^w = 0 $) to an assumption that no feed-through uncertainty is multiplicative with control inputs ones ($ D_y^w \Delta_k D_z^u = 0 $). Moreover, it represents a more general uncertainty structure (correlated block diagonal).

Numerical simulations exemplified that exploiting the internal structure of uncertainties may provide a significant reduction in conservativeness while still yielding controllers which guarantee robust stability and performance.

Although the proposed method provides both necessary and sufficient conditions for the synthesis of GCCs with a more general uncertainty structure, it still provides only sufficient conditions for the synthesis of robust static output-feedback controllers. Besides the underlying limitation that Guaranteed Cost Controllers only provide sufficient conditions for the optimal controller robust synthesis due to the simplifying assumption $ P(\Delta) = P $ on the structure of the cost function.

As future works, we intend to include the extension of such method for constrained linear uncertain systems, hybrid systems and its application to mobile robotic system control.

\bibliographystyle{interact}
\bibliography{refs}
\appendix

\section{Proof of Lemma \ref{lem_structured_uncertainty_set_properties} for the Optimal Guaranteed Cost Control}
\label{appendix_lossless_proof}

In this appendix, we prove Lemma \ref{lem_structured_uncertainty_set_properties}. Towards this goal, we first present supporting theory required for the proof.

\begin{lemma}
For all $ i \in [1, n] $, let $ \mathbb{S}_i $ be lossless sets. Then,%
\begin{equation}
\mathbb{H}\left(\text{conv}\left( \underset{i = 1}{\overset{n}{\bigcup}} \mathbb{S}_i \right) \right) = \mathbb{H}\left( \underset{i = 1}{\overset{n}{\bigcup}} \mathbb{S}_i \right) = \underset{i = 1}{\overset{n}{\bigcap}} \mathbb{H} (\mathbb{S}_i)
\end{equation}%
where $\mathbb{H}(\mathbb{S})$ is the set of matrices $H$ that satisfy property \textit{(iii)} of Definition \ref{def_lossless_set}.

\begin{proof}
Due to the associative property of the involved operations, we will restrict the prove, without loss of generality, to the case $ n = 2 $. We will first prove the second equality, then proceed to prove the first. 

Let $H \in \mathbb{H}(\mathbb{S}_0) \cap \mathbb{H}(\mathbb{S}_1)$. Then,%
\begin{equation}
\forall S \in \mathbb{S}_0 \cup \mathbb{S}_1: \tr(S H) \le 0
\end{equation}%
and%
\begin{equation}
H = \underset{i = 0}{\overset{\text{rank}(H)}{\sum}} \lambda_i \lambda_i^*,\; \forall S \in \mathbb{S}_0 \cup \mathbb{S}_1 : \lambda_i^* S \lambda_i \le 0.
\end{equation}%
Therefore, $ H \in \mathbb{H}(\mathbb{S}_0 \cup \mathbb{S}_1) $ and $ \mathbb{H}(\mathbb{S}_0 \cup \mathbb{S}_1) = \mathbb{H}(\mathbb{S}_0) \cap \mathbb{H}(\mathbb{S}_1) $. 

Consider now the following set: %
\begin{equation}
\mathbb{S}_2 = \text{conv}(\mathbb{S}_0 \cup \mathbb{S}_1) = \{ \alpha_0 S_0 + \alpha_1 S_1 \mid S_0 \in \mathbb{S}_0, S_1 \in \mathbb{S}_1, \alpha_0 \ge 0, \alpha_1 \ge 0\},
\end{equation}%
and $ H \in \mathbb{H}(\mathbb{S}_2) $. Then, from property \textit{(iii)}, we obtain%
\begin{equation}
\forall \alpha_0 \ge 0, \alpha_1 \ge 0, S_0 \in \mathbb{S}_0, S_1 \in \mathbb{S}_1, \alpha_0 \tr(S_0 H) + \alpha_1 \tr(S_1 H) \le 0
\label{eq_lemma_1_prop_1}
\end{equation}%
and%
\begin{equation}
H = \underset{i = 0}{\overset{\text{rank}(H)}{\sum}} \lambda_i \lambda_i^*,\; \alpha_0 \lambda_i^* S_0 \lambda_i + \alpha_1 \lambda_i^* S_1 \lambda_i \le 0.
\label{eq_lemma_1_prop_2}
\end{equation}%
If $ H $ is chosen such that $ H \in \mathbb{H}(\mathbb{S}_0 \cup \mathbb{S}_1) $, then \eqref{eq_lemma_1_prop_1} and \eqref{eq_lemma_1_prop_2} are immediately satisfied. However, if we choose $ H $ such that $ H \not\in \mathbb{H}(\mathbb{S}_0 \cup \mathbb{S}_1) $, it is impossible to satisfy \eqref{eq_lemma_1_prop_1} and \eqref{eq_lemma_1_prop_2} for cases $ \alpha_0 = 1, \alpha_1 = 0$ and $ \alpha_0 = 0, \alpha_1 = 1$. Therefore, $ \mathbb{H}(\text{Conv}(\mathbb{S}_0 \cup \mathbb{S}_1)) = \mathbb{H}(\mathbb{S}_0 \cup \mathbb{S}_1) $.
\end{proof}
\label{lem_convex_hull_lossless}
\end{lemma}

\begin{corollary}
The convex hull of unions of lossless sets is also lossless, since it is a convex cone (properties \textit{(i)} and \textit{(ii)} hold) and property \textit{(iii)} holds from Lemma \ref{lem_convex_hull_lossless}.
\label{cor_convex_hull_lossless}
\end{corollary}

\begin{lemma}
Let $ C \in \Re^{r(p+q) \times n} $, and $ Q \in \Re^{r \times r} $ be a positive semi-definite matrix. Then, the set%
\begin{equation}
\mathbb{S} = \left\{ \left. C^T \begin{bmatrix}
Q \otimes I_p & 0\\
0 & - Q \otimes I_q
\end{bmatrix} C \right| \forall Q \succ 0\right\}
\end{equation}%
is lossless.

\begin{proof}(Based on proof of Lemma 3 of \citep{iwasaki2000generalized})
Properties \textit{(i)} and \textit{(ii)} from Definition \ref{def_lossless_set} are trivial. Therefore, we focus on proving property \textit{(iii)} holds.

Let%
\begin{equation}
S = C^T \begin{bmatrix}
Q \otimes I_p & 0\\
0 & - Q \otimes I_q
\end{bmatrix} C
\end{equation}%
for brevity. Then, from the properties of Kronecker product, there exists a self-similar transformation $ T $ such that%
\begin{equation}
S = C^T T^T \begin{bmatrix}
I_p \otimes Q & 0\\
0 & - I_q \otimes Q
\end{bmatrix} T C.
\end{equation}

Let $ H \in \Re^{n \times n} $ be a positive semi-definite matrix such that $ \tr(HS) \le 0 $ for all $ S \in \mathbb{S} $. Then, $ H $ admits a full rank factorization $ H = G G^T $ where $ G \in \Re^{n \times r} $ where $r = \text{rank}(H) $. We define%
\begin{equation}
\begin{bmatrix}
W \\ Z
\end{bmatrix} = T C G,
\end{equation}%
$ W = [W_1^T, W_2^T, \ldots, W_p^T]^T $, and $ Z = [Z_1^T, Z_2^T, \ldots, Z_q^T]^T $. Then,%
\begin{equation}
\begin{aligned}
\tr(HS) & = \tr(G G^T S) = \tr(G^T S G)\\
&= \tr\left( \begin{bmatrix}
W \\ Z 
\end{bmatrix}^T\begin{bmatrix}
I_p \otimes Q & 0\\
0 & - I_q \otimes Q
\end{bmatrix}
\begin{bmatrix}
W \\ Z 
\end{bmatrix} \right)\\
&= \tr\left(\underset{i = 1}{\overset{p}{\sum}} W_i^T Q W_i - \underset{i = 1}{\overset{q}{\sum}} Z_i^T Q Z_i \right) \\
&= \underset{i = 1}{\overset{p}{\sum}} \tr(W_i^T Q W_i) - \underset{i = 1}{\overset{q}{\sum}} \tr(Z_i^T Q Z_i) \\
&= \underset{i = 1}{\overset{p}{\sum}} \tr(W_i W_i^T Q) - \underset{i = 1}{\overset{q}{\sum}} \tr(Z_i Z_i^T Q) \le 0.
\end{aligned}
\label{eq_lossless_proof_1}
\end{equation}%
Notice that \eqref{eq_lossless_proof_1} must hold for all $ Q \succeq 0 $. Then, since $ W^T W = \underset{i = 1}{\overset{p}{\sum}} W_i^T W_i $ and $ Z^T Z = \underset{i = 1}{\overset{p}{\sum}} Z_i^T Z_i $, \eqref{eq_lossless_proof_1} is equivalent to %
\begin{equation}
W^T W - Z^T Z \preceq 0,
\end{equation}%
which also implies that $ W^T (I_p \otimes Q) W - Z^T (I_q \otimes Q) Z \preceq 0 $, since $ I \otimes Q \succeq 0 $.

Assume that \textit{(iii)} doesn't hold, then there exists $ \xi_i \in \Re^n $ such that $ H = \underset{i = 1}{\overset{r}{\sum}} \xi_i \xi_i^T $, $ \tr(HS) \le 0\; \forall S \in \mathbb{S} $, and $ \xi_i^T S \xi_i > 0 $ for some $ i $ and $ S \in \mathbb{S} $.

Let $ \xi_i = G u_i $ where $ u_i \in \Re^{r} $, $ ||u_i||_2 = 1 $ and $ \underset{i = 1}{\overset{r}{\sum}} u_i u_i^T = I $. Then,%
\begin{equation}
\xi_i^T S \xi_i = u_i^T G^T S G u_i = u_i^T \left( W^T (I_p \otimes Q) W - Z^T (I_q \otimes Q) Z \right) u_i > 0
\end{equation}%
and $ \tr(HS) > 0 $, which contradicts the initial assumption. Therefore, \textit{(iii)} holds and $\mathbb{S}$ is lossless.
\end{proof}
\label{lem_sub_set_lossless}
\end{lemma}

We are now ready to prove Lemma \ref{lem_structured_uncertainty_set_properties}.
\begin{proof}
We first prove property \textit{(i)} ($ \mathbb{S}^u \subseteq \mathbb{S} $), then continue to prove property \textit{(ii)} ($ \mathbb{S} $ is lossless). 

\textit{(i)}: Let $ \Lambda_p = \tau I_{n_p} $ and $ \Lambda_q = \tau I_{n_q} $. Then, $ \mathbb{S} $ reduces to $ \mathbb{S}^u $. Therefore, $ \mathbb{S}^u \subseteq \mathbb{S} $. 

\textit{(ii)}: For all $ i \in [1, s] $, let $ \mathbb{S}_i = \{S^s_i(\Lambda_i) \mid \forall \Lambda_i \succeq 0 \} $, where%
\begin{equation}
S^s_i(\Lambda_i) = 
\begin{bmatrix}
\bar{C}_{z,i} & \bar{D}_{z,i}^w\\ 
0 & I_{n_{ri} n_{qi}}
\end{bmatrix}^T
\begin{bmatrix}
\Lambda_i \otimes I_{n_{pi}} & 0\\ 
0 & -\Lambda_i \otimes I_{n_{qi}}
\end{bmatrix}
(\bullet).
\end{equation}%
Then, we conclude that%
\begin{equation}
\mathbb{S} = \text{conv}\left( \underset{i = 1}{\overset{s}{\bigcup}} \mathbb{S}_i \right).
\end{equation}

Based on Lemma \ref{lem_sub_set_lossless} and Corollary \ref{cor_convex_hull_lossless}, $ S_i $ is lossless for all $ i \in [1, s] $ and S is lossless since it is the convex hull of the union of lossless sets.
\end{proof}

\end{document}